\documentclass[11pt]{amsart}
\usepackage{amsmath}
\usepackage{amsfonts}
\usepackage{amssymb}
\usepackage{amsthm}
\usepackage[english]{babel}

\newtheorem{theo}{Theorem}

\newtheorem{prop}[theo]{Proposition}

\newcommand{\E}{\ensuremath{\mathbb {E}}}
\renewcommand{\P}{\ensuremath{\mathbb {P}}}

\newcommand{\N}{{\mathbb N}}
\newcommand{\C}{{\mathbb C}}
\newcommand{\Q}{{\mathbb Q}}

\renewcommand{\S}{{\mathbb S}}

\def\E{\mathbb{E}}
\def\Z{\mathbb{Z}}

\def\R{\mathbb{R}}

\begin{document}

\title[Good sequences with uncoutable spectrum]{Good sequences with uncountable spectrum and singular asymptotic distribution}

\author{Christophe Cuny}
\author{Fran\c cois Parreau}

\curraddr{Univ Brest, UMR CNRS 6205, LMBA, 6 avenue Victor Le Gorgeu, 29238 Brest}
\email{christophe.cuny@univ-brest.fr}

\curraddr{Universite\'e Sorbonne Paris Nord, LAGA, CNRS, UMR 7539, F-93430, Villetaneuse, France}
\email{parreau@math.univ-paris13.fr}



\begin{abstract}
We construct a good sequence with uncountable spectrum. As an application, we answer to a question of Lesigne, Quas Rosenblatt and Wierdl.
\end{abstract}

\maketitle

\medskip

\medskip

\section{Good sequences with uncountable spectrum}

Let $S=(s_n)_{n\ge 1}$ be an increasing sequence of positive integers. 
We say that $S$ is a \emph{good sequence} if the following limit exists for every $\lambda\in \S^1$ ($\S^1=\{z\in \C\, :\, |z|=1\}$), 
\begin{equation}\label{definition}
c(\lambda)=c_S(\lambda):= \lim_{N\to +\infty}\frac1N\sum_{n=1}^N \lambda^{s_n} \, .
\end{equation}
Equivalently, $S$ is good if,  for every $\lambda \in \S^1$, the following limit exists 
\begin{equation}\label{2nd-definition}\lim_{N\to +\infty}\frac1{\pi_S(N)}\sum_{1\le k\le N, \, k\in S}\lambda^k ,
\end{equation}
where $\pi_S(N)=\# \, (S\cap [1,N])$.

\medskip

Good sequences have been studied by many authors.  See for instance
Rosenblatt and Wierdl \cite{RW} (who introduced that notion), Rosenblatt \cite{Rosenblatt}, 
Boshernitzan, Kolesnik, Quas and Wierdl \cite{BKQW}, 
Lema\'nczyk, Lesigne, Parreau, Voln\'y and Wierdl \cite{LLPVW} or 
Cuny, Eisner and Farkas \cite{CEF}.

\medskip

Given a good sequence $S$, we define its spectrum as the set 
\begin{equation}\label{spectrum}
\Lambda_S:= \{\lambda\in \S^1\, :\, c(\lambda)
\neq 0\}\, .
\end{equation}

\medskip

By theorem 2.22 of \cite{RW} (due to Weyl), for any good sequence $S$, 
$\Lambda_S$ has Lebesgue measure 0. If moreover, 
$S$ has positive upper density, i.e.\ satisfies $\limsup_{N\to +\infty} ({\pi_S(N)}/{N})>0$, then $\Lambda_S$ is countable. See Proposition 2.12 and Corollary 2.13 of \cite{CEF} for a proof based on a result of Boshernitzan published in 
\cite{Rosenblatt}. See also \cite{Kahane} for more general results  of that type. 

\medskip

On another hand, up to our knowledge, no good sequence with 
uncountable spectrum is known.  

\medskip

In \cite{CEF}, good sequences have been studied in connection with  Wiener's lemma. In particular they obtained the following results for good sequences, see their Proposition 2.6 and Theorem 2.10. Recall that if $\tau$ is a finite measure on $\S^1$, then $\hat \tau (n)
=\int_{\S^1} \lambda^n d\tau(\lambda)$, for every $n\in \Z$.

\begin{prop}\label{prop-CEF}
Let $S=(s_n)_{n\ge 1}$ be a good sequence. Then, for every 
probability measure $\mu$ on $\S^1$, we have 
\begin{equation*}\label{wiener}
\frac1N\sum_{n=1}^N |\hat \mu (s_n)|^2 \underset{N\to +\infty}
\longrightarrow \int_{(\S^1)^2}c(\lambda_1\bar \lambda_2)
d\mu(\lambda_1)d\mu(\lambda_2)\, .
\end{equation*}
In particular, if $S$ has countable spectrum and $\mu$ is continuous,
\begin{equation}\label{continuous-measure}
\frac1N\sum_{n=1}^N |\hat \mu (s_n)|^2 \underset{N\to +\infty}
\longrightarrow 0\, .
\end{equation}
\end{prop}
\noindent {\bf Remark.}
\eqref{continuous-measure} implies that $\hat \mu(s_n)$ converges in density to 0, by the Koopman-von Neumann Lemma (see e.g.\ Lemma 2.1 of \cite{CEF}). 

\medskip

The above considerations yield and put into perspective the following 
question: does there exist a good sequence with uncountable 
spectrum ? 

\medskip

We answer positively to that question below. To state the result, we need some more notation.

\medskip

Let $(m_j)_{j\ge 1}$ be an increasing sequence of positive integers, 
such that $m_{j+1}/m_j\ge 3$ for every $j\ge 1$.

\medskip 

We associate with $(m_j)_{j\ge 1}$ the sequence $S=(s_n)_{n\ge 1}$ made out of the integers 
(an empty sum is assumed to be 0)
\begin{equation}\label{mk}\Bigl\{ m_k +\sum_{1\le j\le k-1}\omega_j m_j\, :\, k\ge 1,\, (\omega_1,\ldots , \omega_{k-1})\in \{-1,0,1\}^{k-1}\Bigr\}
\end{equation}
 in increasing order. Notice that our assumption on $(m_j)_{j\ge 1}$ implies that all the integers in \eqref{mk} are positive and distinct.

\medskip 

Denote by $\|\cdot\|$ the distance to the nearest integer: $\|t\|:= \min \{|m-t|\, :\, m\in \Z\}$ for every $t\in \R$.

\medskip

\begin{theo}\label{theo-francois}
Let $(m_j)_{j\ge 1}$ be an increasing sequence of positive integers, 
such that $m_{j+1}/m_j\ge 3$ for every $j\ge 1$ and define $S$ as above. Then, $S$ is good and 
\begin{equation}\label{lambda}
\Lambda:=\bigl\{{\rm e}^{2i\pi \theta}\, :\, \theta \in [0,1)\backslash \Q, \, 
\sum_{j\ge 1} \|m_j\theta\|^2 <\infty\bigr\}\subset \Lambda_S\, .
\end{equation}
\end{theo}

\medskip

\noindent {\bf Proof.} For every $k\ge 1$, consider the following set  of integers
\begin{equation*}
M_k:=\Bigl\{\sum_{1\le j\le k-1}\omega_j m_j\, :\, (\omega_\ell)_{1\le \ell <k}
\in \{-1,0,1\}^{k-1}\Bigr\}\, .
\end{equation*}

For every $k\ge 1$ and every $\theta \in [0,1)$, set 
\begin{align}
\label{Lk-def} L_k(\theta):  &  =\prod_{1\le j\le k-1} \frac13 (1+2\cos(2\pi m_j\theta))
\\
\nonumber&  =\frac1{3^{k-1}}\prod_{1\le j \le k-1}(1+{\rm e}^{-2i\pi  m_j \theta}
+{\rm e}^{2i\pi m_j\theta})\\
\label{Lk}&  =\frac1{3^{k-1}}\sum_{x\in M_k} {\rm e}^{2i \pi x}\, .
\end{align}

Let $\theta\in [0,1)$. As $-1/3\le (1+2\cos(2\pi\theta m_j))/3\le 1$ for all $j$, if $1+2\cos(2\pi\theta m_j)$ is infinitely often non positive, then $(L_k(\theta))_{k\ge 1}$ converges to 0. 

\medskip
Assume now that $1+2\cos(2\pi\theta m_j)>0$ for $j\ge J$, for some integer $J$. Then, the convergence of $(L_k(\theta))_{k\ge 1}$ follows from the convergence of $(\prod_{j=J}^k (1+2\cos(2\pi\theta m_j))/3)_{k\ge J}$ which is clear since we have an infinite product of positive terms less than or equal to 1. Moreover this infinite product converges, i.e.\ the limit is non-zero, if and only if
\[
\sum_{k=J}^{\infty}  \bigl[1-\frac13(1+2\cos(2\pi m_k \theta))\bigr] = \sum_{k=J}^{\infty}\frac23\big(1-\cos (2\pi m_k\theta)\big)<+\infty,
\]
which is equivalent to $\sum_{k=J}^\infty \|m_k\theta\|^2<+\infty$.
\medskip

If $\theta$ is in the set $\Lambda$ defined by \eqref{lambda} the above condition  is satisfied and moreover, as $\theta$ is then irrational, the product $\prod_{j=1}^{J-1} (1+2\cos(2\pi\theta m_j))/3$ does not vanish.

\medskip
Hence in any case $(L_k(\theta))_{k\ge 1}$ converges, say to $L(\theta)$, and $L$ does not vanish on $\Lambda$.

\bigskip

We wish to prove that for every $\theta\in [0,1)$, $(\frac1N\sum_{n=1}^N 
{\rm e}^{2i\pi s_n\theta})_{N\ge 1}$ converges to $L(\theta)$.

\medskip

Let $N\ge 1$. Since $(s_n)_{n\ge 1}$ is the increasing sequence made out of the numbers given by \eqref{mk}, we can write $s_{N+1}= m_{k_N}+\sum_{1\le j\le k_N-1}\omega_j(N) m_j$.

The integers $s_1,\, \ldots , s_N$ may be split into consecutive blocks
\[
m_1+M_1, \ldots, m_{k_N-1}+M_{k_N-1}, \, W_N,
\]
where 
$
W_N=\{\ell \in m_{k_N}+M_{k_N}\, :\, \ell \le s_N\}\, .
$

As each block $M_k$ consists in $3^{k-1}$ integers, we have
\begin{equation}\label{estimerkN}
\frac{3^{k_N-1}-1}{2}\le N< \frac{3^{k_N}-1}{2}
\end{equation}

We may furthermore split $W_N$ into translates of blocks $M_k$. Namely, if $\omega_{k_N-1}(N)\ne -1$, then $W_N$ begins with $m_{k_N}-m_{k_N-1}+M_{k_N-1}$, if $\omega_{k_N-1}(N)= 1$ another block $m_{k_N}+0\times m_{k_N-1}+M_{k_N-1}$ follows, and so on. More precisely, $W_N$ is the disjoint union
$$
W_N=\bigcup_{1\le j \le k_N-1}\, \bigcup_{\omega < \omega_j(N)}
\Big(m_{k_N}+\sum_{\ell= j+1}^{k_N-1}\omega_{\ell}(N)m_{\ell}+\omega m_j+  M_j\Big)\, .
$$

Hence, by \eqref{Lk},
\begin{equation}\label{sum}
\sum_{n=1}^N{\rm e}^{2i\pi s_n\theta}= 
\sum_{j=1}^{ k_N-1} 3^{j-1} {\rm e}^{2i\pi m_j\theta}L_j(\theta)
+ \sum_{j=1}^{k_N-1}\sum_{\omega < \omega_{j}(N)} 3^{j-1}{\rm e}^{2i\pi u_j(\omega)\theta}L_{j}(\theta)\, ,
\end{equation}
where $u_j(\omega)=m_{k_N} +\sum_{\ell= j+1}^{k_N-1}\omega_{\ell}(N)m_\ell+\omega m_j$.

\medskip

Let us first assume that $L(\theta)= 0$. Then we have  
$$
\frac1N \Big|\sum_{n=1}^N{\rm e}^{2i\pi s_n\theta}\Big|\le 
\frac1N \sum_{j=1}^{ k_N-1} 3^{j-1} |L_j(\theta)|
+ \frac1N \sum_{ j=1}^{ k_N-1}\sum_{\omega < \omega_j(N)} 
3^{j -1}|L_{j}(\theta)|\underset{N\to +\infty}
\longrightarrow 0\, ,
$$
where the convergence follows from a Ces\`aro type argument.

\medskip

Assume now that $L(\theta)\neq 0$. Then, ${\rm e}^{2i\pi m_n\theta}\underset{n\to +\infty}\longrightarrow 1$. 

\medskip

Fix $\varepsilon>0$. Let $r\in \N$ be such that ${\rm e}^{-r}<\varepsilon$, and let $d\in \N$ be such that
$|1-{\rm e}^{2i\pi m_j\theta }|<\varepsilon/(r+1)$ and $|L(\theta)-L_j(\theta)|<\varepsilon$ for every $j\ge d$.

\medskip

For every $N$ such that $k_N\ge d+r$, we have on one hand, since $(L_n(\theta))_{n\ge 1}$ is bounded by 1,
\begin{equation}\begin{aligned}\label{premiers_termes}
\sum_{j=1}^{k_N-r-1} 3^{j-1} |{\rm e}^{2i\pi m_j\theta}L_j(\theta)-L(\theta)|
&+ \sum_{j=1}^{k_N-r-1}\sum_{\omega<\omega_j(N)} 3^{j-1}|{\rm e}^{2i\pi u_j(\omega)\theta}L_j(\theta)-L(\theta)|\\
&\le \sum_{j=1}^{k_N-r-1} 3^{j-1} [2+2\times 2]\le 3^{k_N-r}<3^{k_N}\varepsilon.
\end{aligned}\end{equation}
And on the other hand, as $k_N-r\ge d$, when $k_N-r\le j\le k_N$ we have $|1-{\rm e}^{2i\pi m_j\theta }|<\varepsilon/(r+1)$ and
\[
|1-{\rm e}^{2i \pi u_j(\omega)\theta}| \le \sum_{\ell=k_N-r}^{k_N}|1-{\rm e}^{2i\pi m_{\ell}\theta }|< \varepsilon,
\]
for every choice of $\omega$. So,
\begin{equation}\begin{aligned}\label{derniers_termes}
\sum_{j=k_N-r}^{k_N-1} 3^{j-1} |{\rm e}^{2i\pi m_j\theta}L_j(\theta)-L(\theta)|
&+ \sum_{j=k_N-r}^{k_N-1}\sum_{\omega<\omega_j(N)} 3^{j-1}|{\rm e}^{2i\pi u_j(\omega)\theta}L_j(\theta)-L(\theta)|\\
&< \sum_{j=k_N-r}^{k_N-1} 3^{j-1} [2\varepsilon+2\times2\varepsilon]<3^{k_N}\varepsilon.
\end{aligned}\end{equation}
Gathering \eqref{premiers_termes} and \eqref{derniers_termes}, we get
\[
\Big|\sum_{n=1}^N{\rm e}^{2i\pi s_n\theta}- NL(\theta)\Big|<2\cdot 3^{k_N}\varepsilon.
\]
Finally, in view of \eqref{estimerkN},
$$
\limsup_{N\to +\infty} \Big|\frac1N\sum_{n=1}^N{\rm e}^{2i\pi s_n\theta}- L(\theta)\Big|\le 12 \varepsilon
\, ,
$$
and the announced result follows since $\varepsilon$ is arbitrarily small.  \hfill $\square$

\bigskip

It follows from Theorem \ref{theo-francois} that, in order to produce 
a good sequence with uncountable spectrum, it is sufficient to exhibit an increasing sequence of integers $(m_j)_{j\ge 1}$, with $m_{j+1}/m_j\ge 3$ for every $j\ge 1$, and such that the subgroup of $\S^1$
\begin{equation}\label{H2}
H_2=H_2((m_j)_{j\ge 1}):=\{{\rm e}^{2i\pi \theta}\, :\, \theta \in [0,1), \, 
\sum_{j\ge 1} \|m_j\theta\|^2 <\infty\}\, ,
\end{equation}
be uncountable. 

\medskip

It turns out that those type of subgroups have been studied in 
 \cite{HMP} (and \cite{Parreau}).

\medskip

A similar subgroup, defined by $H_1 :=\{{\rm e}^{2i\pi \theta}\, :\, \theta \in [0,1), \, \sum_{j\ge 1} \|m_j\theta\| <\infty\}$, studied in \cite{HMP} in connection with $H_2$, has also been considered by Erd\H{o}s and Taylor \cite{ET} and by Bergelson et al. 
\cite{BJLR} in connection with {IP}-rigidity. 

\medskip

In the above papers, sufficient conditions have been obtained for $H_2$ or $H_1$ to be uncountable.

\medskip

To state the results concerning $H_2$ subgroups, we shall need a strengthening of the lacunarity condition. We say 
that $(m_j)_{j\ge 1}$ satisfies assumption $(A)$ if   one of the conditions $(A_1)$ or $(A_2)$ below is satisfied:

\begin{align*}
&(A_1) \qquad \sum_{j\ge 1}\Big(\frac{m_j}{m_{j+1}}\Big)^2<\infty\,\\
&(A_2)\qquad \forall j\ge 1\quad  m_j|m_{j+1}\quad \mbox{and} \quad  
m_{j+1}/m_j\underset{j\to +\infty}\longrightarrow \infty\, .
\end{align*}

\begin{prop}\label{prop-uncountable-0}
Let $(m_j)_{j\ge 1}$ be a sequence of integers satisfying assumption $(A)$. Then, $H_2((m_j)_{j\ge 1})$ is uncountable.
\end{prop}

The proposition was proved by the second author \cite{Parreau} (see also section 4.2 of \cite{HMP}) under $(A_1)$ (notice that the condition $\inf_{j\in \N}m_{j+1}/m_j\ge 3$ used in \cite{Parreau} and \cite{HMP} is not restrictive for the uncountability of $H_2$). Actually, it is proved in \cite{Parreau} and \cite{HMP} that $H_2$ supports a continuous (singular) probability measure given by a symmetric Riesz product. A proof of the uncountability of $H_2$ under $(A_1)$ can also be derived from the proof of Theorem 5 in \cite{ET}, which states that $H_1$ is uncountable when $\sum_{j\ge 1}m_{j+1}/m_j <\infty$. 

\medskip

Under condition $(A_2)$, the proposition follows from Theorem 3 in \cite{ET} which states that $H_1\subset H_2$ is uncountable. We use their argument below in the proofs of Proposition \ref{measure_on_lambda} and Theorem \ref{theo-christophe}.

\medskip

We are now able to state our main result, which follows in a straightforward way from Proposition \ref{prop-uncountable-0} and Theorem \ref{theo-francois}.

\begin{theo}\label{main-theo}
Let $(m_j)_{j\ge 1}$ be an increasing sequence of positive integers, 
such that $m_{j+1}/m_j\ge 3$ for every $j\ge 1$ and define $S$ as above. If assumption $(A)$ is satisfied then $S$ is a good sequence and it  has uncountable spectrum.
\end{theo}

\medskip

We also derive the following proposition which complements Proposition \ref{prop-CEF}. It can be shown as an abstract consequence of the existence of a good sequence with uncountable spectrum, but we shall give explicit examples.

\begin{prop}\label{measure_on_lambda}
There exist a good sequence $(s_n)_{n\ge 1}$ and a continuous measure 
$\mu$ on $\S^1$, such that $(\frac1N\sum_{n=1}^N |\hat \mu(s_n)|^2)_{N\ge 1}$ converges to some positive number.
\end{prop}

\noindent {\bf Proof.} We construct such a measure for each sequence $S$ associated with a sequence $(m_j)_{j\ge1}$ satisfying ($A_2$) and $\inf_{j\in \N}m_{j+1}/m_j\ge 3$.

Under this assumption, choose a subsequence $(m_{j_k})_{k\ge 1}$ such that $j_1>1$ and $m_{j}/m_{j-1}>2^{k+2}$ for all $j\ge j_k$. For every sequence $\eta=(\eta_k)_{k\ge 1}\in \{0,1\}^{\N^\ast}$, let
\begin{equation*}
\theta(\eta)=\sum_{k= 1}^{\infty}\frac{\eta_k}{m_{j_k}}\;.
\end{equation*}

Given $j\ge1$, let $k$ be the smallest integer such that $j_k> j$. Since $m_j/m_{j_\ell}$ is an integer when $\ell<k$, we have
\begin{equation}\label{m-theta-eta}
\|m_j\theta(\eta)\|\le m_j \sum_{\ell\ge k}\frac{1}{m_{j_\ell}}\le 2\,\frac{m_{j\,}}{m_{j_k}},
\end{equation}
and in particular $\|m_j\theta(\eta)\|\le 1/4$, which yields that all the terms in the products \eqref{Lk-def} are positive.

We also have $\sum_{j<j_k}m_j^2<2 m_{j_k-1}^2$, so if we sum up the $\|m_j\theta(\eta)\|^2$ by blocks from $j_{k-1}$ to $j_{k}-1$ (or from $1$ to $j_1-1$ for the first one), we get that each partial sum is less than $8 (m_{j_k-1}/m_{j_k})^2$,
\[
\sum_{j=1}^{\infty}\|m_j\theta(\eta)\|^2< 8 \sum_{k=1}^{\infty}\frac{m_{j_k-1}^2}{m_{j_k}^2}<\sum_{k=1}^{\infty}\frac{1}{4^{k}}<+\infty.
\]
and $L(\theta(\eta))>0$ follows.

Now, let $\xi =(\xi_j)_{j\ge 1}$ be a sequence of i.i.d.\ random variables with $\P(\xi_1=0)=\P(\xi_1=1)=\frac12$ and let $\mu$ be the probability distribution of ${\rm e}^{2i\pi\theta(\xi)}$. Then, as the mapping $\eta\mapsto {\rm e}^{2i\pi\theta(\eta)}$ is one-to-one, $\mu$ is a continuous probability measure concentrated on $\Lambda_S$ and
\[
\frac1N\sum_{n=1}^N \hat \mu(s_n)=\E\Bigl(\frac1N\sum_{n=1}^N {\rm e}^{2i\pi\theta(\xi)}\Bigr) \to \E\bigl(L(\theta(\xi)\bigr)>0
\quad\mbox{as }N\to+\infty.
\]
Finally, Proposition \ref{prop-CEF} ensures the convergence of $\frac1N\sum_{n=1}^N |\hat \mu(s_n)|^2$ and the positivity of the limit follows the inequality
$\frac1N\sum_{n=1}^N |\hat \mu(s_n)|^2\ge \left|\frac1N\sum_{n=1}^N \hat \mu(s_n)\right|^2$. \hfill $\square$

\medskip

\noindent {\bf Remark.} Under assumption $(A_1)$ and 
$\inf_{j\in \N}m_{j+1}/m_j\ge 3$, the result holds for the measure $\mu$ constructed in \cite{Parreau} or \cite{HMP}. Indeed then $\mu$ is a generalized Riesz product, weak*-limit of products of trigonometric polynomials $P_j$ with coefficients in blocks $\{km_j\,;\,-k_j\le k\le k_j\}$ and $\hat P_j(m_j)=\hat P_j(-m_j)=\cos(\pi/(m_j+2))$. Then for every $s=\sum_{1\le j\le n}\omega_j m_j$ where $|\omega_j|\le k_j$ for all $j$, we have $\hat\mu(s)=\Pi_{1\le j\le n}\,\hat P_j(\omega_j m_j)$ (see \cite{HMP}). From there, the convergence of $\frac1N\sum_{n=1}^N |\hat \mu(s_n)|^2$ and the positivity of the limit can be proven as in Theorem \ref{theo-francois} (we skip the details).

\goodbreak
\section{singular asymptotic distribution}

We now turn to a matter adressed by Lesigne, Quas, Rosenblatt and Wierdl in the preprint \cite{LQRW}. 

\medskip

Let $S=(s_n)_{n\in \N}$ be a good sequence. Let $\lambda\in \S^1$. 
Since $S$ is good, the sequence  $\bigl(\frac1N\sum_{n=1}^N\hat\delta_{\lambda^{s_n}}(m)\bigr)_{N\in\N}=\bigl(\frac1N\sum_{n=1}^N\lambda^{ms_n}\bigr)_{N\in\N}$ converges towards $c(\lambda^m)$ for any integer $m$, that is for any character on $\S^1$, so that $(\frac1N\sum_{n=1}^N\delta_{\lambda^{s_n}})_{N\in\N}$ converges weakly to some probability $\nu_{S,\lambda}$. 

\medskip

Given a probability measure $\nu$ on $\S^1$, if there exists a good sequence $S$ and $\lambda\in \S^1$ such that $\nu_{S, \lambda}=\nu$, one says that $S$ represents the measure $\nu$ at the point $\lambda$.

 \medskip
Lesigne et al. proved several interesting results concerning the  measures that can be represented by a good sequence at some point $\lambda\in \S^1$. 
For instance, they proved that if $\lambda$ is not a root of unity then $\nu_{S,\lambda}$ is continuous (see their Theorem 8.1). They also proved that if a given probability measure $\nu$ on $\S^1$ is not Rajchman (i.e.\ its Fourier coefficients do not vanish at infinity) then, for almost every $\lambda$ with respect to the Haar measure, there does not exist any good sequence representing $\nu$ at $\lambda$ (see their Theorem 8.2). On the opposite, if $\nu$ is absolutely continuous with respect to the Haar measure, then for every $\lambda\in \S^1$ which is not a root of unity there exists a good sequence $S$ representing $\nu$ at $\lambda$ (see their Theorem 9.1).

 \medskip

 The above results raise the following questions. Does there exist a continuous but singular probability measure $\nu$ on $\S^1$  that 
  can be represented by a good sequence? If so, can one take $\nu$ to be non Rajchman?  
 
 \medskip

 It turns out that Theorem \ref{theo-francois} allows to exhibit a good sequence $S$ and a point $\lambda$ such that $\nu_{S,\lambda}$ is a non Rajchman probability measure.

 \medskip
 
\begin{theo}\label{theo-christophe}
 Let $(m_j)_{j\in \N}$ be an increasing sequence of integers satisfying $(A_2)$ and\\ $\inf_{j\in \N}m_{j+1}/m_j\ge 3$, and let $S$ be the sequence associated with it. There are uncountably many $\lambda\in \Lambda_S$ such that the 
 weak*-limit $\nu_{S,\lambda}$ of 
 $(\frac1N\sum_{n=1}^N \delta_{\lambda^{s_n}})_{N\in \N}$ satisfies  $\limsup_{j\to +\infty}|\hat\nu_{S,\lambda}(m_j)|= 1$. 
 \end{theo}
\medskip
\noindent {\bf Proof.} We proceed as in the proof of Proposition \ref{measure_on_lambda}, except that we require a stronger condition on the subsequence $(m_{j_k})_{k\ge 1}$, namely $m_{j}/m_{j-1}>2^{k+2}m_{j_{k-1}}$ for all $j\ge j_k$ if $k>1$.

\medskip
For $\eta\in\{0,1\}^{\N^{\ast}}$, we still define $\theta(\eta)=\sum_{k\ge1}\eta_k/m_{j_k}$.
By the proof of Proposition \ref{measure_on_lambda}, this yields an uncountable family of $\lambda={\rm e}^{2i\pi\theta(\eta)}$ in $\Lambda_S$.

For each such $\theta=\theta(\eta)$ we have $\hat\nu_{S,\lambda}(m)=c({\rm e}^{2i\pi m\theta(\eta)})=L(m\theta)$ for all $m\in\Z$. So, it will be sufficient to show that $L(m_{j_n}\theta)\to 1$ as $n\to+\infty$. Clearly, from the expression of $L(\theta)$ as an infinite product, it is equivalent to prove that
$\sum_{j\ge 1}^{\infty}\|m_{j_n}m_j\theta\|^2$ converges to $0$ as $n\to +\infty$.

\medskip Fix $n>1$. We may apply the inequality (\ref{m-theta-eta}) either to $\|m_j\theta\|$ or to $\|m_{j_n}\theta\|$. For $j<j_n$ we get $\|m_{j_n}m_j\theta\|\le m_j\|m_{j_n}\theta\| \le 2\,m_{j_n}m_j/{m_{j_{n+1}}}$, and in the opposite case $\|m_{j_n}m_j\theta\|\le m_{j_n}\|m_j\theta\|\le 2\,m_{j_n}m_j/{m_{j_k}}$ where $k$ is the smallest integer such that $j_k > j$.
So,
\[
\sum_{j=1}^{j_n-1}\|m_{j_n}m_j\theta\|^2\le 4\, \frac{m_{j_n}^2}{m_{j_{n+1}}^2}\sum_{j=1}^{j_n-1}m_j^2
\le 8\, \frac{m_{j_n}^2}{m_{j_{n+1}}^2}\, m_{j_n-1}^2 < \frac{1}{4^n}\,.
\]
For $j\ge j_n$, we sum again by blocks from $j_{k-1}$ to $j_{k}-1$, for $k>n$,
\[
\sum_{j_{k-1}}^{j_k-1}\|m_{j_n}m_j\theta\|^2
\le 4\, \frac{m_{j_n}^2}{m_{j_{k}}^2}\sum_{j_{k-1}}^{j_k-1}m_j^2
\le 8\, \frac{m_{j_n}^2}{m_{j_{k}}^2}\, m_{j_k-1}^2< \frac{1}{4^k}\frac{m_{j_n}^2}{m_{j_{k-1}}^2} \le \frac{1}{4^k}\,,
\]
and finally
\[
\sum_{j=1}^{\infty}\|m_{j_n}m_j\theta\|^2< \sum_{k=n}^{\infty}\frac{1}{4^k}\to 0\quad\mbox{ as }n\to+\infty.
\]
\hfill $\square$

Let $S$ and $\lambda\in \S^1$ be as in Theorem \ref{theo-christophe} and write $\nu=\nu_{S,\lambda}$. 

\medskip
The property $\limsup_{j\to +\infty}|\hat\nu(m_j)|= 1$ means precisely that $\nu$ is a \emph{Dirichlet 
measure}, see \cite{HMP} and \cite{HMP2} for properties of Dirichlet measures.

In particular there is then a subsequence $(n_j)_{j\ge 1}$ such that $\lambda^{n_j}$ converges towards a constant of modulus 1 in the $L^1(\nu)$ topology, and it follows that any measure absolutely continuous with respect to $\nu$ is itself a Dirichlet measure. 

\medskip
On the other hand, any probability measure absolutely continuous with respect to some Rajchman measure is itself a Rajchman measure. 

\medskip

Hence, we infer that $\nu$ is singular with respect to any Rajchman probability measure 
on $\S^1$.

\medskip

\noindent {\bf Question.}  In view of Theorem \ref{theo-christophe}, one may wonder if it is possible to find a good sequence $S$ and $\lambda\in \S^1$ such that
$$0<\limsup_{n\to +\infty}|\hat \nu_{S,\lambda}(n)|<1\, .$$
Another question is whether one can have $\nu_{S,\lambda}$ Rajchman and singular with respect to the Lebesgue measure.

\medskip

\noindent {\bf Acknowledgement.} We would like to thank Emmanuel Lesigne for  interesting discussions on the topic and a careful reading of a former version of the paper.

\end{document}